\documentclass[11pt]{amsart}
\usepackage{amssymb}

\newcommand{\limfrom}{{\displaystyle\lim_{\longleftarrow}}}    
\newcommand{\limto}{{\displaystyle\lim_{\longrightarrow}}}
\newcommand{\rightlim}{\mathop{\limto}}

\newcommand{\leftlim}{\mathop{\displaystyle\lim_{\longleftarrow}}}
\newcommand{\limfromn}{\leftlim\limits_{\raise3pt\hbox{$n$}}}
\newcommand{\limton}{\rightlim\limits_{\raise3pt\hbox{$n$}}}

\newcommand{\rightlimit}[1]{\mathop{\lim\limits_{\longrightarrow}}\limits%
                   _{\raise3pt\hbox{$\scriptstyle #1$}}}
  
\newcommand{\leftlimit}[1]{\mathop{\lim\limits_{\longleftarrow}}\limits%
                   _{\raise3pt\hbox{$\scriptstyle #1$}}}

\newcommand{\cA}{{\mathcal A}}
\newcommand{\cB}{{\mathcal B}}
\newcommand{\cC}{{\mathcal C}}
\newcommand{\cD}{{\mathcal D}}
\newcommand{\cE}{{\mathcal E}}

\newcommand{\cL}{{\mathcal L}}

\newcommand{\cO}{{\mathcal O}}

\newcommand{\BA}{{\mathbb A}}
\newcommand{\BC}{{\mathbb C}}
\newcommand{\BF}{{\mathbb F}}

\newcommand{\BN}{{\mathbb N}}

\newcommand{\BQ}{{\mathbb Q}}
\newcommand{\BZ}{{\mathbb Z}}

\newcommand{\fm}{{\mathfrak m}}

\newcommand{\iso}{\buildrel{\sim}\over{\longrightarrow}}

\newcommand{\mono}{\hookrightarrow}

\newcommand{\Ltensor}{\buildrel L\over\otimes}

\DeclareMathOperator{\Aut}{{Aut}}
\DeclareMathOperator{\Bad}{{Bad}}
\DeclareMathOperator{\Char}{{char}}

\DeclareMathOperator{\End}{{End}}
\DeclareMathOperator{\et}{{et}}
\DeclareMathOperator{\Ext}{{Ext}}
\DeclareMathOperator{\Fix}{{Fix}}
\DeclareMathOperator{\Fr}{{Fr}}
\DeclareMathOperator{\Gal}{{Gal}}
\DeclareMathOperator{\gr}{{gr}}

\DeclareMathOperator{\Hom}{{Hom}}

\DeclareMathOperator{\IP}{{Irr}}

\DeclareMathOperator{\Ker}{{Ker}}
\DeclareMathOperator{\id}{{id}}
\DeclareMathOperator{\im}{{Im}}

\DeclareMathOperator{\Mor}{{Mor}}

\DeclareMathOperator{\prf}{{prf}}
\DeclareMathOperator{\Pic}{{Pic}}

\DeclareMathOperator{\RGamma}{{R\Gamma}}
\DeclareMathOperator{\RHHom}{{R\bf Hom}}

\DeclareMathOperator{\Spec}{{Spec}}
\DeclareMathOperator{\Supp}{{Supp}}
\DeclareMathOperator{\Top}{{top}}

\theoremstyle{definition}

\numberwithin{equation}{section}

\begin{document}

\thanks{Partially supported by NSF grant DMS-0100108}

\title{On a conjecture of Kashiwara}

\author{Vladimir Drinfeld}

\address{Dept. of Math., Univ. of Chicago, 5734 University 
Ave., Chicago, IL 60637}

\email{drinfeld@math.uchicago.edu}

\dedicatory{To Yu.I.Manin on his 65th birthday}

\maketitle

\section{Introduction}

\subsection{}  \label{gen}

There is a number of deep results on irreducible perverse
sheaves of geometric origin on varieties over $\BC$ proved by
reduction modulo $p$ and a purity argument (see, e.g.,
\S 6.2 of\cite{BBD}). We give a general method to deduce them
{\it without the geometric origin assumption\,} from a
plausible conjecture by  A.J.~de~Jong on local systems
modulo $l$ on varieties over a finite field of
characteristic $p\ne l$. To demonstrate the method we deduce
from de Jong's conjecture the following surprising
conjecture formulated by Kashiwara \cite{K}.

\medskip

\subsection{Conjecture Kash$_{\Top}(\BC )$.}
Suppose that ``algebraic variety'' means ``algebraic variety
over $\BC$" and ``perverse sheaf on $X$'' means ``perverse
sheaf of $F$-vector spaces on $X(\BC )$ equipped with the
usual topology, which is constructible in the Zariski
sense'', where $F$ is a field of characteristic $0$.

1. Let $\pi:X\to Y$ be a proper morphism of algebraic
varieties and $M$ a semisimple perverse sheaf on $X$. Then
the complex $\pi_*(M)$ is ``semisimple'', i.e., isomorphic
to a direct sum of complexes of the form $N_k[k]$ where the
$N_k$'s are semisimple perverse sheaves.

2. In the above situation the hard Lefschetz theorem holds,
i.e., if $u\in H^2(X,\BZ )$ is the class of a relatively
ample line bundle on $X$ then multiplication by $u^k$,
$k>0$, induces an isomorphism
$H^{-k}\pi_*(M)\to H^k\pi_*(M)$.

3. Let $M$ be a semisimple perverse sheaf on an algebraic
variety $X$. Let $f$ be a regular function on $X$. Denote
by $\Psi_f$ the corresponding nearby cycle functor.  Let
$W$ denote the monodromy filtration of $\Psi_f(M)$ (i.e., the
unique exhaustive increasing filtration with the following
property: if $T_u\in\Aut\Psi_f(M)$ is the unipotent
automorphism from the Jordan decomposition of the monodromy
$T\in\Aut\Psi_f(M)$ and $N:=T_u-\id$ then
$NW_k\Psi_f(M)\subset W_{k-2}\Psi_f(M)$ and the morphisms
$N^k:\gr^W_k(\Psi_f(M))\to\gr^W_{-k}(\Psi_f(M))$ are
isomorphisms). Then the perverse sheaf $\gr^W(\Psi_f(M))$ is
semisimple. 

\medskip

{\bf Remark.} Statement 3 implies the
semisimplicity of $\gr^W(\Phi_f(M))$, where $\Phi_f$ is the
vanishing cycle functor. Indeed, if $M$ is irreducible and
$\Supp M\not\subset\{x\in X|f(x)=0\}$ then
$(\Phi_f(M),N)\simeq \im (T-\id:\Psi_f(M)\to\Psi_f(M))$, so 
$\gr^W(\Phi_f(M))$ considered just as a perverse sheaf (with
forgotten grading) is isomorphic to a direct summand of
$\gr^W(\Psi_f(M))$.

\subsection{}
Beilinson, Bernstein, Deligne, and Gabber proved that the
above conjecture holds if $M$ ``has geometric origin'' (see
\cite{BBD}, \S\S6.2.4 -- 6.2.10 for the proof of statements
1 and 2 in this case). We will show that it holds if the rank
of each irreducible component of $M$ on its smoothness locus
is $\le 2$, and in the general case it would follow from a
plausible conjecture formulated by de Jong \cite{dJ}. We
need the following particular case of de Jong's conjecture.

\medskip

\noindent {\bf Conjecture dJ(n).} Let $X$ be an absolutely
irreducible  normal scheme over a finite filed $\BF_q$. Let
$\BF$ be a finite field such that $\Char\BF\ne\Char\BF_q$.
Suppose that
$\rho_0:\pi_1 (X\otimes\bar\BF_q)\to GL(n,\BF )$ is an
absolutely irreducible representation and 
$\rho_t:\pi_1 (X\otimes\bar\BF_q)\to GL(n,\BF [[t]])$ is a
deformation of $\rho_0$ (i.e., a continuous morphism such
that the composition of $\rho_t$ and the evaluation morphism
$GL(n,\BF [[t]])\to GL(n,\BF )$ equals $\rho_0$). If $\rho_t$
extends to a morphism $\pi_1 (X)\to GL(n,\BF [[t]])$ then
the deformation $\rho_t$ is trivial (i.e., 
$\rho_t(\gamma )=g_t\rho (\gamma )g_t^{-1}$, $g_t\in\Ker
(GL(n,\BF [[t]])\to GL(n,\BF ))$, $g_0=1$).

\subsection{Main Theorem} \label{main}
{\it Conjecture dJ($n$) implies Kashiwara's
conjecture for irreducible perverse sheaves whose rank over
the smoothness locus equals $n$.}

\subsection{}
De Jong \cite{dJ} has proved dJ(n) (and, in fact, a
stronger conjecture) for $n\le 2$. His method is to reduce
the statement to the case that $X$ is a smooth projective
curve and then to prove a version of the Langlands
conjecture (more precisely, he proves that given a continuous
representation $\sigma:\pi_1 (X)\to GL(n,\BF ((t))\,)$, 
$n\le 2$, whose restriction to $\pi_1 (X\otimes\bar\BF_q)$
is absolutely irreducible there exists a nonzero 
$\BF ((t))$-valued unramified cusp form on $GL(n)$ over the
adeles of $X$ which is an eigenfunction of the Hecke
operators with eigenvalues related in the usual way with
the eigenvalues of $\sigma (Fr_v)$, $v\in X$). Then he uses
the fact that the space of cusp forms with given central
character is finite-dimensional and therefore the
eigenvalues of Hecke operators belong to $\bar\BF$ if the
central character is defined over $\bar\BF$.

\medskip

\subsection{Remark} The above Conjecture 
Kash$_{\Top}(\BC )$ implies a similar statement for holonomic
$\cD$-modules with regular singularities. Kashiwara \cite{K}
conjectured that the regular singularity assumption is, in
fact, unnecessary (in this case statement 3 is more
complicated, see \cite{K}). I cannot prove this.

\subsection{}
For an algebraically closed field $k$ denote by Kash$_l(k)$
the analog of Conjecture Kash$_{\Top}(\BC )$ with $\BC$
replaced by $k$, the usual topology replaced by the
etale topology and $F$ assumed to be a finite extension of
$\BQ_l$, $l\ne\Char k$. One has
Kash$_{\Top}(\BC )\Rightarrow$Kash$_l(\BC )$; on the other
hand, if Kash$_l(\BC )$ is true for infinitely many primes
$l$ then Kash$_{\Top}(\BC )$ holds (the easy Lemma
\ref{constrlemma} below shows that the field $F$ from
Kash$_{\Top}(\BC )$ can be assumed to be a number field
without loss of generality; then use arguments from
\cite{BBD}, \S6.1.1). One also has 
Kash$_l(\BC )\Rightarrow$Kash$_l(k)$ for every $k$ of
characteristic 0 (indeed, a perverse
$\BQ_l$-sheaf on a variety over $k$ is defined over a
countable subfield
$k_0\subset k$, which can be embedded into $\BC$). In fact,
a specialization argument in the spirit of \cite{BBD},
\S6.1.6 shows that 
Kash$_l(k)\Leftrightarrow$Kash$_l(\bar\BQ )$ for any $k$ of
characteristic 0.

\medskip

\noindent {\bf Remark.} A similar specialization argument
shows that Kash$_l(k)$ is equivalent to Kash$_l(\bar\BF_p)$
for any algebraically closed field $k$ of characteristic
$p>0$. But I cannot prove Kash$_l(\bar\BF_p)$.

\subsection{}
To prove that de Jong's conjecture implies 
Kash$_{\Top}(\BC )$ we use the fact that Kash$_l(\bar\BF_p)$
is true if $M$ comes from a perverse sheaf $M_0$ on a variety
$X_0$ over $\BF_{p^m}$ with
$X_0\otimes_{\BF_{p^m}}\bar\BF_p=X$. This is well known if
$M$ is pure, but according to Lafforgue (\cite{L}, Corollary
VII.8) every absolutely irreducible l-adic perverse sheaf on
$X_0$ becomes pure after tensoring it by a rank 1 sheaf on
$\Spec\BF_{p^m}$ (this follows from the Langlands conjecture
for $GL(n)$ over a functional field proved in \cite{L}).
We also use the moduli of local systems on a complex
variety to get rid of certain singularities (see
\ref{innocent}). This is why I cannot deduce
Kash$_l(\bar\BF_p)$ from de Jong's conjecture even though we
use characteristic $p$ arguments to deduce  
Kash$_{\Top}(\BC )$.

Even if it is possible to prove de Jong's conjecture, it
would be great if somebody finds a direct proof of 
Kash$_{\Top}(\BC )$ or its $\cD$-module version. This has
already been done in some particular cases. In the case that
$Y$ is a point, $X$ is smooth, and $M$ is a local system
statement 2 was proved by C.~Simpson \cite{Si}. Using
$\cD$-modules and Simpson's idea of mixed twistor 
structures C.~Sabbah \cite{Sa} has recently proved statements
2,3 of Kashiwara's conjecture under the assumption that $X$
is smooth and $M$ is a local system extendable to a local
system on a compactification on $X$; he has also proved
statement 1 under the additional assumption that $Y$ is
projective.

\subsection{}
I thank A.~Beilinson for drawing my attention to
Kashiwara's conjectures and numerous discussions. I also
thank A.J. de Jong, who had informed me about his work
\cite{dJ} long before it was published.

\section{Outline of the proof} \label{outline}

In this section we prove Theorem \ref{main} modulo some
lemmas. Their proofs (which are standard) can be found in
\S\ref{constr} and \S\ref{rest}.

\subsection{} \label{prelim}
We can assume that the perverse sheaf $M$ from Conjecture
Kash$_{\Top}(\BC )$ is irreducible and is not supported on
a closed subvariety of $X$ different from $X$. Then there is
a smooth open $j:U\mono X$, $U\ne\emptyset$, such that 
$M=j_{!*}M_U$ for some irreducible lisse perverse sheaf $M_U$
on $U$. Denote by $n$ the rank of $M_U$.

A lisse perverse sheaf on $U(\BC )$ is the same
as a local system on $U$ or a representation of
$\pi_1(U,u)$, $u\in U(\BC )$. So irreducible lisse perverse
sheaves on $U(\BC )$ form an algebraic stack of finite type
over $\BZ$. The corresponding coarse moduli scheme $\IP^U_n$
is of finite type over $\BZ$. The functor
$A\mapsto\IP^U_n(A)$ on the category of commutative rings is
the sheaf (for the etale topology) associated to the
presheaf $\underline{\IP}^U_n$ where $\underline{\IP}^U_n(A)$
is the set of isomorphism classes of rank $n$ locally free
sheaves of $A$-modules $N$ on $U(\BC )$ such that
$N\otimes_Ak$ is irreducible for every field $k$ equipped
with a homomorphism $A\to k$. It is easy to see that if
$\Pic A=0$ and every Azumaya algebra over $A$ is isomorphic
to $\End A^n$ then the map $\underline{\IP}^U_n\to\IP^U_n$ is
bijective. In particular, this is true if $A$ is a complete
local ring with finite residue field.

Notice that locally constant sheaves with finite fibers on 
$U(\BC )$ and $U_{\et}$ are the same. Moreover, if $U$ has a
model $U_E$ over a subfield $E\subset\BC$ (i.e.,
$U=U_E\otimes_E\BC$) then they are the same as locally
constant  sheaves with finite fibers on $(U_E\otimes_E\bar
E)_{\et}$, where $\bar E$ is the algebraic closure of $E$ in
$\BC$, so that $\Gal(\bar E/E)$ acts on the disjoint
union of the completions $(\widehat{\IP}^U_n)_z$ of
$\IP^U_n$ at all possible closed points $z\in\IP^U_n$. 

\subsection{} \label{prelim2}
Fix a closed point
$z\in\IP^U_n$ and a model $U_E$ of $U$ over a finitely 
generated subfield $E\subset\BC$ so that $z$ is 
$\Gal(\bar E/E)$-invariant. Then $\Gal(\bar E/E)$ acts on
$(\widehat{\IP}^U_n)_z$. In \S\ref{rest} we will prove the
following lemma.

\subsection{Lemma} \label{action} 
{\it (a) The action of 
$\Gal(\bar E/E)$ on $(\widehat{\IP}^U_n)_z$
is unramified almost everywhere, i.e., it factors through\,
$\pi_1(\Spec R,\Spec\bar E)$ for some finitely generated
ring $R$ with field of fractions $E$. For such $R$ we have
a well-defined Frobenius automorphism\, $\Fr_v$ of
$(\widehat{\IP}^U_n)_z$ corresponding to a closed point
$v\in\Spec R$ and an embedding of the henselization $R_v$
into $\bar E$: namely $\Fr_v\in\pi_1(\Spec R,\Spec\bar E)$ is
the image of the Frobenius element of $\pi_1(\Spec R,\bar v)$
under the isomorphism 
$\pi_1(\Spec R,\bar v)\iso\pi_1(\Spec R,\Spec\bar E)$
corresponding to the embedding $R_v\mono\bar E\,$.

(b) Assume that Conjecture dJ($n$) holds. Then there exists
$R$ as in (a) such that for every closed point $v\in\Spec R$
and $k\in\BN$ the fixed point scheme 
$\Fix (\Fr_v^k,(\widehat{\IP}^U_n)_z)$ is finite over
$\BZ_l$, where $l$ is the characteristic of the residue
field of $z$.}

\subsection{} \label{Bad}
Now let us prove Theorem \ref{main}.
If $M$ is as in Conjecture Kash$_{\Top}(\BC )$ and $F'$ is a
field containing $F$ then each of the 3 statements of the
conjectures holds for $M\otimes_FF'$ if and only if it holds
for $M$. So for $i=1,2,3$ there is a subset
$\Bad^{\BQ}_i\subset \IP^U_n\otimes\BQ$ such that for every
field $F\supset\BQ$ and every rank $n$ absolutely irreducible
lisse perverse sheaf of $F$-vector spaces $M_U$ on $U$ the
$i$-th statement of Conjecture Kash$_{\Top}(\BC )$ does not
hold for $M=j_{!*}M_U$ if and only if
$M_U\in\Bad_i^{\BQ}(F)$. In  \S\ref{constr} we will prove
the following statement.

\subsection{Lemma} \label{constrlemma}
{\it $\Bad_i^{\BQ}\subset\IP^U_n\otimes\BQ$ is constructible
for $i=1,2,3$.}

\subsection{} \label{innocent}
Denote by $\Bad_i$ the closure of $\Bad_i^{\BQ}$ in
$\IP^U_n$. We have to prove that $\Bad_i^{\BQ}=\emptyset$. If
$\Bad_i^{\BQ}\ne\emptyset$
then there is an open
$V\subset{\IP}^U_n$ such that 
$\Bad_i^{\BQ}\cap (V\otimes\BQ)$ is closed in $V\otimes\BQ$
and $V\cap\Bad_i\subset V$
is non-empty and smooth over $\BZ$ (if equipped with
the reduced scheme structure). Take a closed point 
$z\in V\cap\Bad_i$ and denote by $(\widehat{\Bad}_i)_z$ the
preimage of $\Bad_i$ in $(\widehat{\IP}^U_n)_z$. Fix $E$ and
$U_E$ as in \ref{prelim2}.

\subsection{Lemma} \label{invariance}
{\it $(\widehat{\Bad}_i)_z$ is $\Gal (\bar E/E)$-stable if
$E$ is big enough in the following sense. We assume that
$j:U\mono X$ comes from a morphism $j_E:U_E\mono X_E$
of schemes over $E$. We also assume that in the cases
$i=1,2\;$ $\pi :X\to Y$ comes from a morphism
$\pi_E:X_E\to Y_E$ of schemes over $E$ and in the case $i=3$
the function $f$ is defined over $E$. Besides, in the case
$i=2$ we assume that our relatively ample line bundle on $X$
comes from a bundle $\cL_E$ on $X_E$. }

\begin{proof}
Let us consider $(\widehat{\IP}^U_n)_z$ as a scheme rather
than a formal scheme. It suffices to show that
if $F$ is a finite extension of $\BQ_l$ and $O\subset F$ is
the ring of integers then the subset 
$T\subset\Mor (\Spec O,(\widehat{\IP}^U_n)_z)$ consisting of
morphisms $f:\Spec O\to (\widehat{\IP}^U_n)_z$ such that
$f(\Spec F)\subset (\widehat{\Bad}_i)_z$ is 
$\Gal (\bar E/E)$-invariant. A morphism 
$f:\Spec O\to (\widehat{\IP}^U_n)_z$ defines an
isomorphism class of lisse perverse sheaves of
$O$-modules $M_U$ on $U(\BC )$, and 
$f(\Spec F)\subset (\widehat{\Bad}_i)_z$ if and only if the 
$i$-th statement of Conjecture Kash$_{\Top}(\BC )$ does not
hold for the irreducible perverse sheaf
$M=j_{!*}(M_U\otimes_OF)$. We can consider $M_U$ (resp.
$M$) as a perverse $O$-sheaf (resp. a perverse
$F$-sheaf) on the scheme $U_E\otimes_E\bar E$ (resp.
$X_E\otimes_E\bar E$) rather than on the topological space
$U(\BC )$ (resp. on $X(\BC )$). So $T$ is
$\Gal (\bar E/E)$-invariant.
\end{proof}

\subsection{Lemma} \label{existence}
{\it Let $E$ be as in \ref{invariance} and $R$, $l$ as in
\ref{action}(b). Let $v$ be a closed point of $\Spec R$ and
$i\in\{ 1,2,3\}$. Then the fixed point scheme
$\Fix_k:=\Fix (\Fr_v^k,(\widehat{\Bad}_i)_z)$ is finite and
flat over $\BZ_l$ for every $k\in\BN$. For some $k\in\BN$ it
is not empty.}

\begin{proof}
We chose $R$ so that finiteness is clear. We chose $z$
so that $(\widehat{\Bad}_i)_z$ is smooth over $\BZ_l$,
therefore $\Fix_k\subset (\widehat{\Bad}_i)_z$ is defined by
$d$ equations, where $d$ is the dimension of
$(\widehat{\Bad}_i)_z$ over $\BZ_l$. So finiteness
implies that $\Fix_k$ is a complete intersection over $\BZ_l$
and therefore flat over $\BZ_l$ (this argument was used,
e.g., in \S3.14 of \cite{dJ}). For some
$k\in\BN$ the automorphism $\Fr_v^k$ acts identically on the
residue field of $z$. Then $\Fix_k$ is not empty. 
\end{proof}

The above lemma contradicts the following
statement, which will be proved in~\S\ref{rest}.

\subsection{Lemma} \label{nopoints}
{\it Let $E$ be as in \ref{invariance} and $R$ as in
\ref{action}(a). For each $i\in\{ 1,2,3\}$ there is a
non-empty open
$W\subset\Spec R$ such that for every closed point
$v\in W$ and every $k\in\BN$ the fixed point scheme 
$\Fix (\Fr_v^k,(\widehat{\Bad}_i)_z)$ has no
$\bar\BZ_l$-points, where $\bar\BZ_l$ is the ring of
integers of $\bar\BQ_l$.}

\section{Constructibility} \label{constr}

In this section all rings are assumed to be
commutative and Noetherian. Lemma \ref{constrlemma} can be
reformulated as follows.

\subsection{Lemma} \label{constrlemma2}
{\it Let $N$ be a locally free sheaf of modules on $U(\BC )$
over an integral $\BQ$-algebra $A$ such that $N\otimes_Ak$ is
irreducible for all fields $k$ equipped with a homomorphism
$A\to k$. Let $M_U:=N[d]$, $d:=\dim U$, be the corresponding
perverse sheaf. Denote by $K$ the field of fractions of $A$.
Let $i\in\{ 1,2,3\}$. 

(a) If the $i$-th statement of Conjecture
Kash$_{\Top}(\BC )$ holds for the perverse sheaf 
$M=j_{!*}M_U\otimes_AK$ then there is an $f\in
A\setminus\{ 0\}$ such that it holds for
$j_{!*}(M_U\otimes_Ak)$ for 
every field $k$ equipped
with a homomorphism $A_f\to k$. 

(b) If the $i$-th statement of Conjecture
Kash$_{\Top}(\BC )$ does not hold for
$M=j_{!*}M_U\otimes_AK$ then there is an 
$f\in A\setminus\{ 0\}$ such that for every field $k$ 
equipped with a homomorphism $A_f\to k$ it does not hold for
$j_{!*}(M_U\otimes_Ak)$.}

\medskip

The proof will be given in \ref{prove} after some preparatory
lemmas.

\subsection{} 
Let $X$ be an algebraic variety over $\BC$. As usual,
$D(X(\BC ),A)$ denotes the derived category of sheaves of
$A$-modules on $X(\BC )$ and 
$D^b_c(X(\BC ),A)\subset D(X(\BC ),A)$ is the full
subcategory of complexes with constructible cohomology
sheaves (the notion of constructible set being understood in
the sense of algebraic geometry). Recall that the perverse
t-structure on $D^b_c(X(\BC ),A)$ is defined as follows:
$\null^pD^{\ge 0}_c(X(\BC ),A)$ (resp. 
$\null^pD^{\le 0}_c(X(\BC ),A)$) is the full subcategory of
$D^b_c(X(\BC ),A)$ consisting of complexes $C$ such that 
every irreducible subvariety 
$Y\buildrel{\nu}\over{\hookrightarrow}X$ has a
non-empty open subset $Y'\subset Y$ such that
$H^i\nu^!C|_{Y'}=0$ (resp.$H^i\nu^*C|_{Y'}=0$) for 
$i<2\dim S$ (resp. for $i>-2\dim S$). We will write
$\null^pD^{\ge 0}$, $\null^pD^{\le 0}$ instead of
$\null^pD^{\ge 0}_c$, $\null^pD^{\le 0}_c$. Let 
$\null^pD_{\ge 0}(X(\BC ),A)\subset D^b_c(X(\BC ),A)$ denote
the full subcategory of complexes $C$ such that
$C\Ltensor_A N$ belongs to $\null^pD^{\ge 0}(X(\BC ),A)$ for
all finitely generated $A$-modules $N$.

\subsection{Lemma} \label{strat}
{\it Suppose we have a stratification 
$X=\bigcup_{\nu} X_{\nu}$. Denote by $i_{\nu}$ the embedding
$X_{\nu}\mono X$.

(i) $C\in\null^pD^{\ge 0}(X(\BC ),A)$ if and only if
$i_{\nu}^!C\in\null^pD^{\ge 0}(X_{\nu}(\BC ),A)$ for all
$\nu$.

(ii) $C\in\null^pD_{\ge 0}(X(\BC ),A)$ if and only if
$i_{\nu}^!C\in\null^pD_{\ge 0}(X_{\nu}(\BC ),A)$ for all
$\nu$.}

\begin{proof}
(i) is obvious. (ii) follows from (i) because the functor
$\Ltensor_A N$ commutes with $i_{\nu}^!$
\end{proof}

\subsection{Lemma} \label{basechange}
{\it If $C\in\null^pD_{\ge 0}(X(\BC ),A)$ then 
$A'\Ltensor_A C\in\null^pD_{\ge 0}(X(\BC ),A)$ for every
$A$-algebra $A'$ (not necessarily finite over $A$).}

\begin{proof}
Using \ref{strat}(ii) we reduce the proof to the case
where $X$ is smooth and $C$ is lisse. This case is equivalent
to that of $X=\Spec\BC$, which is well known and easy.
\end{proof}

\subsection{}
Let $M$ be a perverse sheaf of $A$-modules. If $N$ is a
finitely generated $A$-module the perverse sheaf
$N\otimes_AM$ is defined by 
$N\otimes_AM:=\tau_{\ge 0}N\Ltensor_AM$, where $\tau_{\ge 0}$
denotes perverse truncation. Quite similarly one defines
$A'\otimes_AM$ for an $A$-algebra $A'$. 

A perverse sheaf of $A$-modules $M$ is said to be {\it
flat\,} if $M\in\null^pD_{\ge 0}(X(\BC ),A)$; this is
equivalent to exactness of the functor 
$N\mapsto N\otimes_AM$ on the category of finitely generated
$A$-modules $N$.

\subsection{Lemma} \label{generic}
{\it Let $A$ be an integral ring.

(i) If $M$ is a perverse sheaf of $A$-modules
then there is an $f\in A\setminus\{ 0\}$ such that
$A_f\otimes_AM$ is flat over $A_f$.

(ii) If $C\in\null^pD^{\ge 0}(X(\BC ),A)$ then there is an 
$f\in A\setminus\{ 0\}$ such that
$A_f\otimes_AC\in\null^pD_{\ge 0}(X(\BC ),A)$.}

\begin{proof}
Clearly (ii)$\Rightarrow$(i). The proof of (ii) is 
similar to that of \ref{basechange}.
\end{proof}

\subsection{Lemma} \label{generic2}
{\it Let $A$ be an integral ring and $K$ its field of
fractions. Let $M$ be a perverse sheaf of $A$-modules on
$X(\BC )$. Let $j:U\mono X$ be an open embedding. Suppose
that $M\otimes_AK=j_{!*}j^*(M\otimes_AK)$. Then there is an
$f\in A\setminus\{ 0\}$ such that for every field $k$
equipped with a homomorphism $A_f\to k$ one has
$M\otimes_Ak=j_{!*}j^*(M\otimes_Ak)$.}

\begin{proof}
Consider the embedding $i:X\setminus U\to X$. The equality
$M\otimes_AK=j_{!*}j^*(M\otimes_AK)$ means that
$i^*M\otimes_AK\in\null^pD^{<0}(X(\BC ),K)$ and
$i^!M\otimes_AK\in\null^pD^{>0}(X(\BC ),K)$. Localizing $A$
one can assume that $i^*M\in\null^pD^{<0}(X(\BC ),A)$ and
$i^!M\in\null^pD^{>0}(X(\BC ),A)$. Moreover, by 
\ref{generic} we can assume that $M$ is flat and
$i^!M\in\null^pD_{>0}(X(\BC ),A)$. Then for every field $k$
equipped with a homomorphism $A\to k$ one has
$M\otimes_Ak=M\Ltensor_Ak$, so 
$i^*(M\otimes_Ak)=(i^*M)\Ltensor_Ak\in
\null^pD^{<0}(X(\BC ),k)$ and 
$i^!(M\otimes_Ak)=i^!(M\Ltensor_Ak)=(i^!M)\Ltensor_Ak$ is in
$\null^pD^{>0}(X(\BC ),k)$ by  \ref{basechange}.
Therefore $M\otimes_Ak=j_{!*}j^*(M\otimes_Ak)$.
\end{proof}

\subsection{Lemma} \label{generic4}
{\it Let $A,K,M,X$ be as in \ref{generic2}. Suppose that
$\Char K=0$.

(i) If the perverse sheaf $M\otimes_AK$ is semisimple then
there is an $f\in A\setminus\{ 0\}$ such that $M\otimes_Ak$
is semisimple for every morphism from $A_f$ to a field $k$.

(ii) If $M\otimes_AK$ is not semisimple
then there is an $f\in A\setminus\{ 0\}$ such that
$M\otimes_Ak$ is not semisimple for every morphism from
$A_f$ to a field $k$.}

\begin{proof}
(i) We can assume that $M\otimes_AK$ is irreducible.
Localizing $A$ we can assume that the support of $M$ equals
the support of $M\otimes_AK$. So we can assume that this
support equals $X$. Then there is a smooth open subset
$j:U\mono X$ such that $M\otimes_AK=j_{!*}j^*(M\otimes_AK)$
and $j^*M\otimes_AK$ is lisse and irreducible. Then
$j^*M\otimes_A\bar K$ is semisimple (here we use that 
$\Char K=0$), so localizing $A$ we can assume that
$j^*M\otimes_Ak$ is semisimple for every morphism from $A$
to a field $k$. It remains to apply \ref{generic2}.

(ii) Let $0\to N_K\to M_K\to P_K\to 0$ be a nontrivial
extension of perverse sheaves with $M_K=M\otimes_AK$. 
Localizing $A$ we can assume that it comes from an
exact sequence $0\to N\to M\to P\to 0$ of flat perverse
sheaves of $A$-modules on $X(\BC )$. It defines an element 
$u\in H^1(C)$, $C:=R\Hom (P,N)$, such that the image of $u$
in $H^1(C)\otimes_AK$ is nonzero. For every morphism from
$A$ to a field $k$ the sequence of perverse sheaves
$0\to N\otimes_Ak\to M\otimes_Ak\to P\otimes_Ak\to 0$ is
still exact and its class in
$\Ext^1(P\otimes_Ak,N\otimes_Ak)=H^1(C\Ltensor_Ak)$ is the
image of $u$. As $C$ is quasi-isomorphic to a finite complex
of finitely generated $A$-modules, there exists 
$f\in A\setminus 0$ such that for every morphism from $A_f$
to a field $k$ the image of $u$ in $H^1(C\Ltensor_Ak)$ is
nonzero.
\end{proof}

\subsection{Lemma} \label{generic5}
{\it Let $A,K$ be as in  \ref{generic4} and
$\cC\in D^b_c(X(\BC ),A)$.

(i) Suppose that $\cC\otimes_AK$ is ``semisimple'', i.e.,
isomorphic to a direct sum of complexes of the form $N_j[j]$,
where the $N_j$'s are semisimple perverse sheaves. Then there
is an $f\in A\setminus\{ 0\}$ such that $\cC\Ltensor_Ak$ is
``semisimple'' for every morphism from $A_f$ to a field $k$.

(ii) If $\cC\otimes_AK$ is not ``semisimple'' then there is
an $f\in A\setminus\{ 0\}$ such that $\cC\Ltensor_Ak$ is not
``semisimple'' for every morphism from $A_f$ to a field $k$.}

\begin{proof}
(i) follows from  \ref{generic4}(i). Let us prove (ii).
By  \ref{generic}(i) we can assume that the perverse
cohomology sheaves $H^i\cC$ are flat over $A$. Then
$H^j(\cC\Ltensor_Ak)=H^j\cC\otimes_Ak$, so if
$H^j\cC\otimes_AK$ is not semisimple for some $j$ then by
 \ref{generic4} there is an
$f\in A\setminus\{ 0\}$ such that $H^j(\cC\Ltensor_Ak)$ is
not semisimple for every morphism from $A_f$ to a field $k$.
If $H^j\cC\otimes_AK$ is semisimple for all $j$ then there
exists $j$ such that the canonical element
$u\in\Hom (\tau_{>j}\cC ,(\tau_{\le j}\cC )[1])$ 
has nonzero image in 
$\Hom (\tau_{>j}\cC ,(\tau_{\le j}\cC )[1])\otimes_AK$.
Just as in the proof of  \ref{generic4}(ii) this
implies that there is an $f\in A\setminus\{ 0\}$ such that
the image of $u$ in $\Hom (\tau_{>j}\cC\Ltensor_Ak,
(\tau_{\le j}\cC\Ltensor_Ak)[1])$ is nonzero for every
morphism from $A_f$ to a field $k$.
\end{proof}

\subsection{Lemma} \label{generic6}
{\it Let $A,K$ be as in  \ref{generic2} and
$\varphi: \cC_1\to\cC_2$ a morphism in $D^b_c(X(\BC ),A)$.
If the morphism
$\varphi_K:\cC_1\otimes_AK\to\cC_2\otimes_AK$ induced by
$\varphi$ is an isomorphism then there is an 
$f\in A\setminus\{ 0\}$ such that
$\varphi_k: \cC_1\Ltensor_Ak\to\cC_2\Ltensor_Ak$
is an isomorphism for every morphism from $A_f$ to a
field $k$. If $\varphi_K$ is not an isomorphism then there
is an $f\in A\setminus\{ 0\}$ such that $\varphi_k$ is not an
isomorphism for every morphism from $A_f$ to a field $k$.}
\hfill \qedsymbol

\subsection{Proof of Lemma \ref{constrlemma2}} \label{prove}

By \ref{generic}(i) and
\ref{generic2}, localizing $A$ we can assume that $j_{!*}M_U$
is $A$-flat and $(j_{!*}M_U)\otimes_Ak=j_{!*}(M_U\otimes_Ak)$
for every morphism from $A$ to a field $k$. 

To prove Lemma
\ref{constrlemma2} for $i=1$ apply \ref{generic5} to
$\cC =\pi_*j_{!*}M_U$ and notice that
$\pi_*j_{!*}(M_U\otimes_Ak)=(\pi_*j_{!*}M_U)\Ltensor_Ak$.
To prove the lemma for $i=2$ notice that by \ref{generic}(i)
one can assume that the perverse sheaves $H^r\pi_*j_{!*}M_U$,
$r\in\BZ$, are $A$-flat and so
$H^r\pi_*j_{!*}(M_U\otimes_Ak)=
(H^r\pi_*j_{!*}M_U)\otimes_Ak$.
Then apply \ref{generic6}.

Now let us prove the lemma for $i=3$. As $\Psi_f$ is exact,
$L:=\Psi_fj_{!*}M_U$ is an $A$-flat perverse sheaf and
$\Psi_fj_{!*}(M_U\otimes_Ak)=L\otimes_Ak$. As
$\dim\End (L\otimes_AK)<\infty$ and $\Char K=0$, the
monodromy $T\in\Aut(L\otimes_AK)$ satisfies $p(T)^m=0$ for
some $m\in\BN$ and some monic polynomial $p$ with nonzero
discriminant $D$ and constant term $c$. Localizing $A$ we can
assume that $D,c$ are invertible in $A$ and the equality
$p(T)^m=0$ holds in $\End L$. There is a unique
$q\in A[t]/(p(t))$ such that $\,q\equiv 1\mod p(t)$ and
$p(t/q)=0$. Define $T_u\in\Aut L$ by $T_u:=q(T)$, then $T_u$
becomes equal to the unipotent part of $T$ after any base
change $A\to k$, where $k$ is a field. There is a unique
exhaustive increasing filtration $W$ on $L$ such that
$(T-1)W_kL\subset W_{k-2}L$ and the morphisms
$(T-1)^k:\gr^W_kL\to\gr^W_{-k}L$ are isomorphisms (this holds
for a unipotent automorphism in any abelian category).
By \ref{generic}(i), localizing $A$ we can assume that
$\gr^WL$ is $A$-flat. Then for every morphism from $A$ to
a field $k$ the filtration on $L\otimes_Ak$ induced by $W$
is the weight filtration and $\gr
(L\otimes_Ak)=(\gr^WL)\otimes_Ak$. It remains to apply
\ref{generic4}. \hfill \qedsymbol

\section{Very good models} \label{good}

In this section we remind the results of \cite{BBD} used
in \S\ref{rest}.

\subsection{} \label{notat}
Let $E$ be a field finitely generated over a
prime field, $X$ a scheme of finite type over $E$ and
$M_1,\ldots,M_k\in D^b_c(X,\BF )$, where $\BF$ is a finite
field whose characteristic $l$ is different from that of
$E$. Let $R\subset E$ be a subring of finite type over $\BZ$
whose field of fractions equals $E$.

A {\it model\,} of $X$ over $R$ is a scheme $X_R$ of finite
type over $R$ with an isomorphism $X_R\otimes_RE\iso X$. A
model of $(X,M_1,\ldots,M_k)$ over $R$ is a collection
$(X_R,(M_1)_R,\ldots,(M_k)_R)$ where $X_R$ is a model of $X$
over $R$ and $(M_i)_R$ is an object of $D^b_c(X_R,\BF )$
whose pull-back to $X$ is identified with $M_i$.

For every $R$ and $E$ as above every collection
$(X,M_1,\ldots,M_k)$ has a model over $R$. Given a model
$(X_R,(M_1)_R,\ldots,(M_k)_R)$ we denote by $X_u$ the fiber 
of $X_R$ over a geometric point $u$ of $\Spec R$ and by
$(M_i)_u$ the $*$-restriction of $M_i$ to $X_u$. We write
$X_{\bar E}$, $(M_i)_{\bar E}$ instead of $X_{\Spec \bar E}$,
$(M_i)_{\Spec \bar E}$ and denote by $p_R$ the morphism
$X_R\to\Spec R$. We need the following ad hoc definition.

\subsection{Definition} \label{goodness}
A model $(X_R,(M_1)_R,\ldots,(M_k)_R)$ is {\it good\,} if

1) $l$ is invertible in $R$;

2) for every geometric point $u$ of $\Spec R$ and every
$i,j$ the morphism
\begin{equation} \label{1}
((p_R)_*\RHHom ((M_i)_R,(M_j)_R))_u\to
 R\Hom ((M_i)_u,(M_j)_u)
\end{equation}
is an isomorphism and the morphism
\begin{equation} \label{2}
((p_R)_*\RHHom ((M_i)_R,(M_j)_R))_u\to
 R\Hom ((M_i)_{\bar E},(M_j)_{\bar E})
\end{equation}
induced by an embedding of the strict henselization $R_u$
into $\bar E$ is also an isomorphism.

\subsection{Remark} \label{Otoo}
In the above definition we understand
$\RHHom$ as $\RHHom_{\BF}$ and $R\Hom:=R\Hom_{\BF}$. But if
$O$ is a local ring with residue field $\BF$ then the
analogs of (\ref{1}) and (\ref{2}) with $\RHHom$ replaced by
$\RHHom_O$ and $R\Hom$ replaced by $R\Hom_O$ are still
isomorphisms (because $\RHHom_O((M_i)_R,(M_j)_R)$ equals
$\RHHom ((M_i)_R,(M_j)_R)\otimes R\Hom_O(\BF ,\BF )$, etc.).
\medskip

\subsection{}
According to \S 6.1 of \cite{BBD}, {\it every model
$(X_R,(M_1)_R,\ldots,(M_k)_R)$ becomes good after a base
change of the form $R\to R_f$, 
$f\in R\setminus\{ 0\}$.} More precisely, after replacing $R$
by some $R_f$ the following properties (which are
stronger than the above property 2) hold:

a) for every geometric point $u$ of $\Spec R$ and every
$i,j$ the morphism 
$(\RHHom ((M_i)_R,(M_j)_R))_u$ 
$\to$ 
$\RHHom ((M_i)_u,(M_j)_u)$ 
is an isomorphism;

b) the complexes $(p_R)_*K_{ij}$, 
$K_{ij}:=\RHHom ((M_i)_R,(M_j)_R)$, are lisse
and for every geometric point $u$ of $\Spec R$ and every
$i$,$j$ the morphism 
$((p_R)_*K_{ij})_u$ 
$\to$ 
$\RGamma (X_u,(K_{ij})_u)$
is an isomorphism.

This is shown in \S 6.1 of \cite{BBD} by reducing to the
case where each $(M_i)_R$ is the extension by zero of a local
system on a locally closed subscheme of $X_R$ and using
Theorem 1.9 from \cite{De3}.

\subsection{}
Suppose we have a model $(X_R,(M_1)_R,\ldots,(M_k)_R)$. Let
$R_u$ denote the strict henselization of $R$ at a geometric
point $u$ of $\Spec R$. Choose an algebraic closure 
$\bar E\supset E$ and an embedding $R_u\mono\bar E$. Let
$O$ be a local Artinian ring whose residue
field $O/\fm $ is a finite extension of $\BF$. Let
$D^{\{M_i\}}(X_u,O)\subset D_c^b(X_u,O)$ be the thick 
triangulated subcategory generated by the
complexes $(M_i)_u\otimes_{\BF}O/\fm$ (according to 
\cite{Ve}, a triangulated subcategory $\cB$ of a triangulated
category $\cA$ is {\it thick\,} if every object of $\cA$
which is a direct summand of an object of
$\cB$ belongs to $\cB$). Let 
$D_{\{M_i\}}(X_u,O)\subset D_c^b(X_u,O)$ denote the
full subcategory of complexes $C$ such that
$C\Ltensor_OO/\fm\in D^{\{M_i\}}(X_u,O/\fm )$. Notice that
$D_{\{M_i\}}(X_u,O)\subset D^{\{M_i\}}(X_u,O)\cap
D_{\prf}(X_u,O)$, where $D_{\prf}(X_u,O)\subset D_c^b(X_u,O)$
is the full subcategory of complexes of finite Tor-dimension.
(In fact, one can prove
\footnote{If 
$C\in D^{\{M_i\}}(X_u,O)$ then for every
$N\in\BN$ one has an exact triangle 
$K\to C\Ltensor_OO/\fm\to K'$ such that
$K\in D^{\{M_i\}}(X_u,O/\fm )$ and $H^jK'=0$ for $j>-N$
(put $K:=C\Ltensor_OP$, where $P$ is a perfect complex of
$O$-modules such that $H^0P=O/\fm$ and $H^jP=0$ if $j\ne 0$
and $j>-N'$). If one also has $C\in D_{\prf}(X_u,O)$ then for
$N$ big enough the morphism $C\Ltensor_OO/\fm\to K'$ is zero,
so $C\Ltensor_OO/\fm$ is a direct summand of 
$K\in D^{\{M_i\}}(X_u,O/\fm )$ and therefore
$C\Ltensor_OO/\fm\in D^{\{M_i\}}(X_u,O/\fm )$.}    
that
$D_{\{M_i\}}(X_u,O)=D^{\{M_i\}}(X_u,O)\cap D_{\prf}(X_u,O)$,
but we do not need this fact). One also has the similar
categories $D^{\{M_i\}}(X_R\otimes_RR_u,O)$,
$D_{\{M_i\}}(X_R\otimes_RR_u,O)$, 
$D^{\{M_i\}}(X_{\bar E},O)$, and $D_{\{M_i\}}(X_{\bar E},O)$.
The following lemma and its proof is a version of \S6.1.9 of
\cite{BBD}.

\subsection{Lemma} \label{goodlemma}
{\it If a model
$(X_R,(M_1)_R,\ldots,(M_k)_R)$ is good then the functors
\begin{equation} \label{equiv1}
D_{\{M_i\}}(X_R\otimes_RR_u,O)\to 
D_{\{M_i\}}(X_{\bar E},O),
\end{equation}
\begin{equation} \label{equiv2}
D_{\{M_i\}}(X_R\otimes_RR_u,O)\to
D_{\{M_i\}}(X_u,O)
\end{equation}
and their analogs for $D^{\{M_i\}}$ are equivalences, so one
gets equivalences of triangulated categories}
\begin{equation}\label{equiv}
D_{\{M_i\}}(X_{\bar E},O)\iso 
D_{\{M_i\}}(X_u,O)
\end{equation}
\begin{equation}\label{equiv3}
D^{\{M_i\}}(X_{\bar E},O)\iso 
D^{\{M_i\}}(X_u,O)
\end{equation}
\begin{proof}
It suffices to prove the lemma for $D^{\{M_i\}}$. To this
end, use Remark \ref{Otoo} and the fact that every idempotent
endomorphism of an object of the derived category of sheaves
comes from a direct sum decomposition (see \cite{Ne}, ch.~1,
Proposition 1.6.8.)
\end{proof}

\subsection{Remark} \label{equivariance}
Consider the point $v\in\Spec R$ corresponding to the
geometric point $u$. The henselization
$R_v$ of $R$ at $v$ is embedded into $\bar E$; denote by
$E_v$ its field of fractions. Then $\Gal (\bar E/E_v)$
acts on $R_u$, $\bar E$, and $u$, so it acts on the
categories $D_{\{M_i\}}(X_R\otimes_RR_u,O)$, 
$D_{\{M_i\}}(X_{\bar E},O)$, and
$D_{\{M_i\}}(X_u,O)$ (an action of a monoid $\Gamma$
on a category $C$ is a monoidal functor from $\Gamma$ to the
monoidal category of functors $C\to C$). The equivalence
(\ref{equiv}) is $\Gal (\bar E/E_v)$-equivariant because it
is a composition of two $\Gal (\bar E/E_v)$-equivariant
equi\-va\-lences.

\subsection{} \label{verygood}
Now suppose that the $M_i$'s are perverse
sheaves whose pull-backs to $X_{\bar E}$ are semisimple.
A good model $(X_R,(M_1)_R,\ldots,(M_k)_R)$ is said to be
{\it very good\,} if $(M_i)_u$ is a semisimple perverse
sheaf for every $i$ and every geometric point $u$ of 
$\Spec R$. In this situation if $(M_i)_{\bar E}$ is
absolutely irreducible (i.e., 
$(M_i)_{\bar E})\otimes_{\BF}\bar\BF$ is irreducible) then
$(M_i)_u$ is absolutely irreducible for every geometric
point $u$ of $\Spec R$ (because (\ref{equiv}) induces an
isomorphism $\End (M_i)_{\bar E}\simeq\End(M_i)_u$).

\medskip

In the case of a very good model $D^{\{M_i\}}(X_u,O)$
consists of all objects of $D_c^b(X_u,O)$ such that all
irreducible components of the reduction modulo $\fm$ of their
perverse cohomology sheaves occur in
$((M_1)_u\oplus\ldots\oplus (M_k)_u)\otimes_{\BF}O/\fm$.
There is a similar description of $D^{\{M_i\}}(X_{\bar
E},O)$. {\it So the equivalences (\ref{equiv3}) and
(\ref{equiv}) corresponding to a very good model send
perverse sheaves to perverse sheaves and the same is true
for the equivalences inverse to (\ref{equiv}) and
(\ref{equiv}).}

Using the principles explained in \S 6.1.7 of \cite{BBD} one
shows that {\it every model $\,(X_R,(M_1)_R,\ldots,(M_k)_R)$
becomes very good  after a base change of the form
$R\to R_f$.}

\subsection{} Now let $O$ be a complete discrete valuation
ring whose residue field $O/\fm $ is a finite extension of
$\BF$. According to Deligne's definition of the l-adic
derived category (see 
\S2.2.14 of \cite{BBD} or \S1.1.2 of \cite{De4}),
$D_c^b(X_u,O)$ is the inverse limit of
$D_{\prf}(X_u,O/\fm^r)$, $r\in\BN$. We define
$D_{\{M_i\}}(X_u,O)$ to be the inverse limit of
$D_{\{M_i\}}(X_u,O/\fm^r )$. Clearly
$D_{\{M_i\}}(X_u,O)\subset D_c^b(X_u,O)$ is the full
subcategory of complexes $C\in D_c^b(X_u,O)$ such that
$C\Ltensor_OO/\fm\in D^{\{M_i\}}(X_u,O/\fm)$. Same is true
for $X_{\bar E}$ and $X_R\otimes_RR_u$. By \ref{goodlemma},
in the case of a good model we have the equivalences
(\ref{equiv1}) -- (\ref{equiv}). It easily follows from
\ref{verygood} that in the case of a very good model the
equivalence (\ref{equiv}) sends perverse sheaves to perverse
sheaves and the same is true for the equivalence inverse to
(\ref{equiv}).

\section{A lemma on nearby cycles}

\subsection{}  \label{notat2}
We keep the notation from \ref{notat}, but now we suppose
that $\Char E=0$. Let $f\in H^0 (X,\cO_{X})$ and 
$M\in D^b_c (X,\BF)$. Let $(X_R,M_R,f_R)$ be a model for
$(X,M,f)$ (i.e., $(X_R,M_R)$ is a model for $(X,M)$ and
$f_R\in H^0 (X_R,\cO_{X_R})$ extends $f$). For each
geometric point $u$ of $\Spec R$ we have the nearby cycle
complex $\Psi_{f_u}(M_u)\in D^b_c (Y_u,\BF)$, where
$f_u=f_R|_{X_u}$ and $Y_u$ is the fiber of the subscheme
$Y_R\subset X_R$ defined by $f_R=0$. Lemma \ref{SGA7} below
essentially says that after localizing $R$ the complexes
$\Psi_{f_u}(M_u)$ come from a single object of 
$D^b_c(X_R,\BF)$. To formulate the lemma precisely we need
some notation. 

Define $X_R[f_R^{1/m}]\subset X_R\times\BA^1$ by the equation
$f_R(x)=t^m$, $x\in X_R$, $t\in\BA^1$. The direct image of
the constant sheaf on $X_R[f_R^{1/m}]$ with fiber $\BF$ with
respect to the projection $X_R[f_R^{1/m}]\to X_R$ is denoted
by $\cE_m$. If $m|m'$ then $\cE_m$ is embedded into
$\cE_{m'}$. Denote by $\cE^{(m)}$ the direct limit of
$\cE_{ml^n}$, $n\in\BN$, where $l:=\Char\BF$. Finally, put
$\Psi_{f_R}^{(m)}(M_R):=
(\nu_R)_*\nu_R^*(M_R\otimes\cE^{(m)})|_{Y_R}$, where
$\nu_R:X_R\setminus Y_R\to X_R$ is the embedding. Clearly
$\Psi_{f_R}^{(ml)}(M_R)=\Psi_{f_R}^{(m)}(M_R)$. For
every geometric point $u$ of $\Spec R$ one has the base
change morphism
\begin{equation} \label{basechange2}
(\Psi_{f_R}^{(m)}(M_R))_u\to  \Psi_{f_u}^{(m)}(M_u)
\end{equation}
and the obvious morphism
\begin{equation} \label{obvious}
\Psi_{f_u}^{(m)}(M_u)\to\Psi_{f_u}(M_u) \; .
\end{equation}
Let $\kappa_u$ be the residue field of $u$ and
$I:=\Gal (\overline{\kappa_u((t))}/\kappa_u((t))\,)$. If
$m,l$ are invertible on $u$ and coprime then
$\Psi_{f_u}^{(m)}(M_u)=\Psi_{f_u}(M_u)^{I_m}$, where $I_m$
is the unique normal subgroup of $I$ such that
$I/I_m\simeq\BZ_l\times\BZ/m\BZ$. The following ad hoc
definition will be used in \ref{i=3}.

\subsection{Definition} \label{adhoc}

A model $(X_R,M_R,f_R)$ of $(X,M,f)$ is {\it $\Psi_f$-good}
if there exists $m\in\BN$ such that

a) $m^{-1},l^{-1}\in R$,

b) the morphisms (\ref{basechange2}) and (\ref{obvious})
are isomorphisms for all $u$,

c) for every $\overline{m}\in m\BN$ the morphism
$\Psi_{f_R}^{(m)}(M_R)\to\Psi_{f_R}^{(\overline{m})}(M_R)$ is
an isomorphism over 
$Y_R\otimes_{\BZ}\BZ [\overline{m}^{-1}]$,

d) the cohomology sheaves of $\Psi_{f_R}^{(m)}(M_R)$ are
constructible.

\subsection{Remarks} \label{rem}
(i) Maybe d) holds automatically.

(ii) If properties a) -- d) hold for $(X_R,M_R,f_R)$ and some
$m$ then for every $m'\in m\BN$ they hold for
$(X_{R'},M_{R'},f_{R'})$ and $m'$, where $R':=R[1/m']$,
$X_{R'}:=X_R\otimes_RR'$, $f_{R'}:=f_R|_{X_{R'}}\,$.

\subsection{Lemma} \label{SGA7}
{\it Every model $(X_R,M_R,f_R)$ of $(X,M,f)$ becomes
$\Psi_f$-good after a base change of the form $R\to R_g$,
$g\in R\setminus \{ 0 \}$.}

\begin{proof}
Using resolution of singularities in characteristic $0$
reduce the proof to the case that 
$X_R=\Spec R[t_1,\ldots t_n]$, 
$f_R=t_1^{k_1}\ldots t_n^{k_n}$, and $M_R$ is the

\noindent constant
sheaf (cf. the proof of Theorem 2.3.1 of \cite{De2}). In this
case take $m$ to be the l.c.m. of $k_1,\ldots ,k_n$ and
perform the base change $R\to R[l^{-1}m^{-1}]$ (cf. the proof
of Theorem 3.3 of \cite{De1}).
\end{proof}

\section{Proof of Lemmas \ref{action} and \ref{nopoints}}
\label{rest}

We fix a model $U_E$ of $U$ over a finitely 
generated subfield $E\subset\BC$. Let $\BF$ denote the
residue field of $z$ (which is finite), and
$l:=\Char\BF$. While proving Lemmas
\ref{action} and \ref{nopoints} we can replace $E$ by its
finite extension, so we can assume that $z$ is the
isomorphism class of a lisse perverse $\BF$-sheaf $M_{U_E}$
on $U_E$. The pull-back $M_{U_{\bar E}}$ of $M_{U_E}$ to 
$U_{\bar E}:=U_E\otimes_E\bar E$ is absolutely irreducible.
We will use the notions of ``good'' and ``very good'' from
\ref{goodness} and \ref{verygood}.

Given a finitely generated ring $R$ with field of fractions
$E$ and a closed point $v\in\Spec R$ we choose an embedding
of the henselization $R_v$ into $\bar E$. We will use this
embedding in the situation of \ref{action}(a). The field of
fractions of $R_v$ is denoted by $E_v$. The embedding
$R_v\to \bar E$ defines a geometric point $\bar v$ of $\Spec
R$ and an embedding of the strict henselization $R_{\bar v}$
into $\bar E$ (if $E_{\bar v}\subset \bar E$ is the maximal
extension of $E_v$ unramified at $v$ and $L$ is the integral
closure of $R_v$ in $E_{\bar v}$ then $\bar v$ is defined to
be the closed point of $L$; then $R_{\bar v}=L$).

\subsection{Proof of Lemma \ref{action}}

Choose a very good model $(U_R,M_{U_R})$ of $(U_E,$
$M_{U_E})$ over a finitely generated ring $R$ with field of
fractions $E$ such that $U_R$ is smooth over $R$ and
$M_{U_R}[-d]$ is a local system, $d:=\dim U$. Then $R$ has
the properties required in Lemma \ref{action}. Indeed, let
$v\in\Spec R$ be a closed point. Our $(\widehat{\IP}^U_n)_z$
is the base of the universal deformation of the local system
$M_{U_{\bar E}}[-d]$. Let $M_{U_{\bar v}}$ be the pull-back
of $M_{U_R}$ to $U_{\bar v}:=U_R\otimes_R\bar v$; as
explained in \ref{verygood}, $M_{U_{\bar v}}$ is  absolutely
irreducible. By \ref{verygood}, for every Artinian local ring
$O$ with residue field $\BF$ the equivalence (\ref{equiv3})
(with $X$ replaced by $U$) induces an equivalence between the
$O$-linear category of perverse sheaves of $O$-modules on
$U_{\bar E}$ and that on $U_{\bar v}$. So we get a canonical
isomorphism between $(\widehat{\IP}^U_n)_z$ and the base of
the universal deformation of $M_{U_{\bar v}}[-d]$. It is
$\Gal (\bar E/E_v)$-equivariant by \ref{equivariance}. So
the action of $\Gal (\bar E/E_v)$ on $(\widehat{\IP}^U_n)_z$
is unramified, and Conjecture dJ(n) implies that for every
$k\in\BN$ the fixed point scheme of
$\Fr_v^k:(\widehat{\IP}^U_n)_z\to (\widehat{\IP}^U_n)_z$ is
finite.

\subsection{Proof of Lemma \ref{nopoints}}
Recall that in the lemma $i$ denotes the number of a 
statement in Kashiwara's conjecture (see \ref{Bad}). We
assume that $E$ is big enough in the sense of
\ref{invariance}. We fix $j_E:U_E\mono X_E$, in the cases
$i=1,2$ we fix $\pi_E:X_E\to Y_E$, and in the case $i=2$ we
also fix $\cL_E$ (see \ref{invariance}).

\subsubsection{} \label{notice}

Let $M_1,\ldots,M_k$ be the irreducible components of the
perverse cohomology sheaf $H^0(j_E)_!M_U$. The pull-backs
of $M_1,\ldots,M_k$ to $X_{\bar E}:=X_E\otimes_E\bar E$ are
semisimple. 

Notice that if $O$ is a complete discrete valuation ring
whose residue field $O/\fm$ contains $\BF$ and $\tilde M$ is
a perverse sheaf of $O$-modules on $U_{\bar E}$ such that 
$\tilde M/\fm \tilde M\simeq M_{U_{\bar E}}
\otimes_{\BF}(O/\fm )$ 
then the perverse sheaf
$(j_{\bar E})_{!*}\tilde M/\fm(j_{\bar E})_{!*}\tilde M$ is a
quotient of $H^0(j_{\bar E})_!M_U\otimes_{\BF}(O/\fm )$ and
so each of its irreducible components occurs in 
$((M_1)_{\bar E}\oplus\ldots (M_k)_{\bar E})
\otimes_{\BF}(O/\fm )$. 

\subsubsection{} 

 Let $R$ be as in \ref{action}(a). Localizing $R$
(i.e., replacing $R$ by $R_f$ for some nonzero $f\in R$) we
can choose a very good model $(X_R,(M_1)_R,\ldots,(M_k)_R)$
of $(X_E,M_1,\ldots,M_k)$. 

\subsubsection{Proof for $i=1,2$} \label{i=12}
Let $N_1,\ldots,N_m$ be the irreducible components of all
perverse cohomology sheaves of the complexes
$(\pi_E)_*M_1,\ldots,(\pi_E)_*M_k$. Localizing $R$ we can
choose a very good model
$(Y_R,(N_1)_R,\ldots,(N_m)_R)$ of $(Y_E,N_1,\ldots,$ $N_m)$.
Further localizing $R$ we can extend 
$\pi_E:X_E\to Y_E$ to a proper morphism $\pi_R :X_R\to Y_R$.
In the case $i=2$ after another localization we can
extend $\cL_E$ to a relatively ample line bundle on $X_R$.
Further localizing $R$ we can assume that 
$(\pi_R)_*(M_1)_R,\ldots ,(\pi_R)_*(M_k)_R$ belong to the
triangulated subcategory of $D^b_c(Y_R,\BF )$ generated by
$(N_1)_R,\ldots,(N_m)_R$. We claim that after these 
localizations the fixed point schemes
$\Fix (\Fr_v^k,(\widehat{\Bad}_i)_z)$, $i\in\{ 1,2\}$, have
no $\bar\BZ_l$-points for every closed point $v\in\Spec R$
and every $k\in\BN$. Let us prove this for $i=1$ (the case
$i=2$ is quite similar).

Suppose there exists a $\Fr_v^k$-invariant
$\bar\BZ_l$-point of $(\widehat{\Bad}_1)_z)$. It comes from
a $\Fr_v^k$-invariant $O$-point $\xi$ of
$(\widehat{\Bad}_1)_z)$, where $O$ is the ring of integers
in some finite extension $F\supset\BQ_l$. We have
$[O/\fm:\BF]<\infty$, where $\fm$ is the maximal ideal of $O$
and $\BF$ is the residue field of $z$. Our $\xi$ corresponds
to an $O$-flat lisse perverse sheaf $\tilde M$ on 
$U_{\bar E}$ such that 
$\tilde M/\fm \tilde M\simeq  M_{U_{\bar E}}
\otimes_{\BF}(O/\fm )$ and the complex
$(\pi_{\bar E})_*(j_{\bar E})_{!*}\tilde M\otimes_OF$ is not
semisimple (here $M_{U_{\bar E}}$ is the pull-back of
$M_{U_E}$ to $U_{\bar E}$ and the
morphisms $j_{\bar E}:U_{\bar E}\mono X_{\bar E}$, 
$\pi_{\bar E}:X_{\bar E}\to Y_{\bar E}$ are induced by
$j_E:U_E\mono X_E$, $\pi_E:X_E\to Y_E$). Besides, 
$(\sigma^k)^*\tilde M\simeq\tilde M$, where 
$\sigma\in\Gal(\bar E/E_v)$ is a preimage of the Frobenius
element of $\pi_1(\Spec R,\bar v)$.

The perverse sheaf
$(j_{\bar E})_{!*}\tilde M\otimes_OF$ is absolutely
irreducible, the complex
$(\pi_{\bar E})_*(j_{\bar E})_{!*}\tilde M\otimes_OF$ is not
semisimple, and 
$(\sigma^k)^*(j_{\bar E})_{!*}\tilde M\simeq
(j_{\bar E})_{!*}\tilde M$. 
By \ref{notice}
$(j_{\bar E})_{!*}\tilde M\in D_{\{M_i\}}(X_{\bar E},O)$, so
$(\pi_{\bar E})_*(j_{\bar E})_{!*}\tilde M \in
D_{\{N_i\}}(Y_{\bar E},O)$. 
We have a commutative diagram of
$\Gal (\bar E/E_v)$-equivariant functors
\begin{equation} \label{diagram1}
\begin{array}{lll}
D_{\{M_i\}}(X_{\bar E},O)&\iso & D_{\{M_i\}}(X_{\bar v},O)\\
\,\downarrow& & \, \downarrow\\
D_{\{N_i\}}(Y_{\bar E},O)&\iso & D_{\{N_i\}}(Y_{\bar v},O)
\end{array}
\end{equation}
where the horizontal arrows are the equivalences
(\ref{equiv}) and the vertical ones are $(\pi_{\bar E} )_*$
and $(\pi_{\bar v} )_*\,$. Let $P\in D_{\{M_i\}}(X_u,O)$ be
the image of $(j_{\bar E})_{!*}\tilde M$. Then $P$ is a
perverse sheaf on $X_{\bar v}$ such that $P\otimes_OF$ is 
absolutely irreducible but the complex 
$(\pi_{\bar v})_*P\otimes_OF$ is not semisimple. Besides,
the isomorphism class of $P$ is invariant with respect to
$\pi_1(v_k,\bar v)$, $v_k:=\Spec \kappa_k$, where $\kappa_k$
is the extension of order $k$ of the residue field of $v$
inside the residue field of $\bar v$. So $P$ is the pull-back
of a perverse sheaf $P_0$ on $X_{v_k}:=X_R\otimes_Rv_k$.
According to Corollary VII.8 of \cite{L} and Corollary 5.3.2
of \cite{BBD}, $P_0\otimes_OF$ becomes pure after tensoring
it by a rank 1 sheaf on $v_k$. So by \S5.1.14 and \S5.4.6 of
\cite{BBD} $(\pi_{\bar v} )_*P\otimes_OF$ is semisimple, and
we get a contradiction.

\subsubsection{Proof for $i=3$}   \label{i=3}  
Let $Y_E\subset X_E$ denote the subscheme $f=0$. Let
$N_1,\ldots,N_r$ be the irreducible components of the
perverse sheaves 
$\Psi_f(M_1)$,\ldots,$\Psi_f(M_k)$.
Localizing $R$ we can assume that $f$ extends to a regular
function $f_R$ on $X_R$.
Further localizing $R$ we can choose a very good model
$(Y_R,(N_1)_R,\ldots,(N_r)_R)$ of $(Y_E,N_1,\ldots,N_r)$
such that $Y_R$ is the closed subscheme of $X_R$ defined by
$f_R=0$. By \ref{SGA7}, after another localization all
the models $(Y_E,N_j,f)$, $1\le j\le r$, become
$\Psi_f$-good in the sense of \ref{adhoc}. By \ref{rem}(ii)
we can assume that the number $m$ from the definition of
``$\Psi_f$-good'' does not depend on $j$.
Further localizing
$R$ we can assume that 
$\Psi_{f_R}^{(m)}(M_1)_R,\ldots ,\Psi_{f_R}^{(m)}(M_k)_R$
are in the triangulated subcategory of $D^b_c(Y_R,$ $\BF )$
generated by
$(N_1)_R,\ldots,(N_r)_R$. 
We claim that after these 
localizations 
$\Fix (\Fr_v^k,(\widehat{\Bad}_3)_z)$ has
no $\bar\BZ_l$-points for every closed point $v\in\Spec R$
and every $k\in\BN$. 

We will explain only the parts of the proof that are not
quite the same as in \ref{i=12}. First, instead of the
semisimplicity theorem from \cite{BBD} we use Gabber's
theorem (\S5.1.2 of \cite{BB}). Second, the role of
(\ref{diagram1}) is played by the diagram 
\begin{equation} \label{diagram2}
\begin{array}{lll}
D_{\{M_i\}}(X_{\bar E},O)&\iso & D_{\{M_i\}}(X_{\bar v},O)\\
\,\downarrow& & \, \downarrow\\
D_{\{N_i\}}(Y_{\bar E},O[[\Gamma]])&\iso 
& D_{\{N_i\}}(Y_{\bar v},O[[\Gamma]])
\end{array}
\end{equation}
Here $O[[\Gamma]]$ is the
completed group algebra of $\Gamma:=\BZ_l(1)$, i.e., the
projective limit of $O[\Gamma/l^n\Gamma]$, $n\in\BN$. We
define the ``$O$-constructible derived category'' 
$D^b_{O\mbox{-}c}(Y_{\bar E},O[[\Gamma]])$ as 
$\limfrom\limto D^b_c(Y_{\bar E},
(O/\fm^s)[\Gamma/l^n\Gamma])$ ($\limfrom$ with respect to
$s$ and $\limto$ with respect to $n$).
We denote by $D_{\{N_i\}}(Y_{\bar E},O[[\Gamma]])$ the full
subcategory of objects 
$C\in D^b_{O\mbox{-}c}(Y_{\bar E},O[[\Gamma]])$ such that the
image of $C$ in $D^b_c(Y_{\bar E},O)$ belongs to
$D_{\{N_i\}}(Y_{\bar E},O)$. The upper horizontal arrow in
(\ref{diagram2}) is the equivalence (\ref{equiv}) and the
lower one is its analog for
$O[[\Gamma]]$. The vertical arrows in (\ref{diagram2}) are
$\Psi_f$ and $\Psi_{f_{\bar v}}\,$. The
$O[[\Gamma]]$-structure on the nearby cycles comes from the
identification of $\Gamma$ with the Sylow $l$-subgroup of
the tame quotient of $\Gal(\overline{K((t))}/K((t))\,)$,
where $K$ is an algebraically closed field, $\Char K\ne l$.
Just as in \ref{i=12}, $\Gal (\bar E/E_v)$ acts on the four
categories from (\ref{diagram2}), and the arrows are
$\Gal (\bar E/E_v)$-equivariant functors. To construct a
canonical isomorphism between the two functors
$D_{\{M_i\}}(X_{\bar E},O)\to 
D_{\{N_i\}}(Y_{\bar v},O[[\Gamma]])$ from (\ref{diagram2})
insert the functor 
$\Psi_{f_{R_{\bar v}}}^{(m)}:
D_{\{M_i\}}(X_R\otimes_RR_{\bar v},O)\to
D_{\{N_i\}}(Y_R\otimes_RR_{\bar v},O[[\Gamma]])$
as a vertical arrow
between the two vertical arrows of (\ref{diagram2}). Here
$\Psi_{f_{R_{\bar v}}}^{(m)}$ is the analog of
$\Psi_{f_R}^{(m)}$ with $R$ and $X_R$ replaced by
$R_{\bar v}$ and $X_R\otimes_RR_{\bar v}$.

\end{document}